
\magnification=1200

\hsize=13.50cm
\vsize=18cm
\parindent=12pt   \parskip=0pt
\pageno=1


\hoffset=15mm    
\voffset=1cm    


\ifnum\mag=\magstep1
\hoffset=-2mm   
\voffset=.8cm   
\fi


\pretolerance=500 \tolerance=1000  \brokenpenalty=5000

\catcode`\@=11

\font\eightrm=cmr8         \font\eighti=cmmi8
\font\eightsy=cmsy8        \font\eightbf=cmbx8
\font\eighttt=cmtt8        \font\eightit=cmti8
\font\eightsl=cmsl8        \font\sixrm=cmr6
\font\sixi=cmmi6           \font\sixsy=cmsy6
\font\sixbf=cmbx6


\font\tengoth=eufm10       \font\tenbboard=msbm10
\font\eightgoth=eufm10 at 8pt      \font\eightbboard=msbm10 at 8 pt
\font\sevengoth=eufm7      \font\sevenbboard=msbm7
\font\sixgoth=eufm7 at 6 pt        \font\fivegoth=eufm5

 \font\tencyr=wncyr10       
\font\eightcyr=wncyr10 at 8 pt      
\font\sevencyr=wncyr10 at 7 pt      
\font\sixcyr=wncyr10 at 6 pt


\skewchar\eighti='177 \skewchar\sixi='177
\skewchar\eightsy='60 \skewchar\sixsy='60


\newfam\gothfam           \newfam\bboardfam
\newfam\cyrfam

\def\tenpoint{%
  \textfont0=\tenrm \scriptfont0=\sevenrm \scriptscriptfont0=\fiverm
  \def\rm{\fam\z@\tenrm}%
  \textfont1=\teni  \scriptfont1=\seveni  \scriptscriptfont1=\fivei
  \def\oldstyle{\fam\@ne\teni}\let\old=\oldstyle
  \textfont2=\tensy \scriptfont2=\sevensy \scriptscriptfont2=\fivesy
  \textfont\gothfam=\tengoth \scriptfont\gothfam=\sevengoth
  \scriptscriptfont\gothfam=\fivegoth
  \def\goth{\fam\gothfam\tengoth}%
  \textfont\bboardfam=\tenbboard \scriptfont\bboardfam=\sevenbboard
  \scriptscriptfont\bboardfam=\sevenbboard
  \def\bb{\fam\bboardfam\tenbboard}%
 \textfont\cyrfam=\tencyr \scriptfont\cyrfam=\sevencyr
  \scriptscriptfont\cyrfam=\sixcyr
  \def\cyr{\fam\cyrfam\tencyr}%
  \textfont\itfam=\tenit
  \def\it{\fam\itfam\tenit}%
  \textfont\slfam=\tensl
  \def\sl{\fam\slfam\tensl}%
  \textfont\bffam=\tenbf \scriptfont\bffam=\sevenbf
  \scriptscriptfont\bffam=\fivebf
  \def\bf{\fam\bffam\tenbf}%
  \textfont\ttfam=\tentt
  \def\tt{\fam\ttfam\tentt}%
  \abovedisplayskip=12pt plus 3pt minus 9pt
  \belowdisplayskip=\abovedisplayskip
  \abovedisplayshortskip=0pt plus 3pt
  \belowdisplayshortskip=4pt plus 3pt
  \smallskipamount=3pt plus 1pt minus 1pt
  \medskipamount=6pt plus 2pt minus 2pt
  \bigskipamount=12pt plus 4pt minus 4pt
  \normalbaselineskip=12pt
  \setbox\strutbox=\hbox{\vrule height8.5pt depth3.5pt width0pt}%
  \let\bigf@nt=\tenrm       \let\smallf@nt=\sevenrm
  \normalbaselines\rm}

\def\eightpoint{%
  \textfont0=\eightrm \scriptfont0=\sixrm \scriptscriptfont0=\fiverm
  \def\rm{\fam\z@\eightrm}%
  \textfont1=\eighti  \scriptfont1=\sixi  \scriptscriptfont1=\fivei
  \def\oldstyle{\fam\@ne\eighti}\let\old=\oldstyle
  \textfont2=\eightsy \scriptfont2=\sixsy \scriptscriptfont2=\fivesy
  \textfont\gothfam=\eightgoth \scriptfont\gothfam=\sixgoth
  \scriptscriptfont\gothfam=\fivegoth
  \def\goth{\fam\gothfam\eightgoth}%
  \textfont\cyrfam=\eightcyr \scriptfont\cyrfam=\sixcyr
  \scriptscriptfont\cyrfam=\sixcyr
  \def\cyr{\fam\cyrfam\eightcyr}%
  \textfont\bboardfam=\eightbboard \scriptfont\bboardfam=\sevenbboard
  \scriptscriptfont\bboardfam=\sevenbboard
  \def\bb{\fam\bboardfam}%
  \textfont\itfam=\eightit
  \def\it{\fam\itfam\eightit}%
  \textfont\slfam=\eightsl
  \def\sl{\fam\slfam\eightsl}%
  \textfont\bffam=\eightbf \scriptfont\bffam=\sixbf
  \scriptscriptfont\bffam=\fivebf
  \def\bf{\fam\bffam\eightbf}%
  \textfont\ttfam=\eighttt
  \def\tt{\fam\ttfam\eighttt}%
  \abovedisplayskip=9pt plus 3pt minus 9pt
  \belowdisplayskip=\abovedisplayskip
  \abovedisplayshortskip=0pt plus 3pt
  \belowdisplayshortskip=3pt plus 3pt
  \smallskipamount=2pt plus 1pt minus 1pt
  \medskipamount=4pt plus 2pt minus 1pt
  \bigskipamount=9pt plus 3pt minus 3pt
  \normalbaselineskip=9pt
  \setbox\strutbox=\hbox{\vrule height7pt depth2pt width0pt}%
  \let\bigf@nt=\eightrm     \let\smallf@nt=\sixrm
  \normalbaselines\rm}

\tenpoint


\def\pc#1{\bigf@nt#1\smallf@nt}         \def\pd#1 {{\pc#1} }


\catcode`\;=\active
\def;{\relax\ifhmode\ifdim\lastskip>\z@\unskip\fi
\kern\fontdimen2  -1.2 \fontdimen3 \string;}

\catcode`\:=\active
\def:{\relax\ifhmode\ifdim\lastskip>\z@\unskip\fi\penalty\@M\ \fi\string:}

\catcode`\!=\active
\def!{\relax\ifhmode\ifdim\lastskip>\z@
\unskip\fi\kern\fontdimen2  -1.1 \fontdimen3 \string!}

\catcode`\?=\active
\def?{\relax\ifhmode\ifdim\lastskip>\z@
\unskip\fi\kern\fontdimen2  -1.1 \fontdimen3 \string?}

\def\^#1{\if#1i{\accent"5E\i}\else{\accent"5E #1}\fi}
\def\"#1{\if#1i{\accent"7F\i}\else{\accent"7F #1}\fi}

\frenchspacing


\newtoks\auteurcourant      \auteurcourant={\hfil}
\newtoks\titrecourant       \titrecourant={\hfil}

\newtoks\hautpagetitre      \hautpagetitre={\hfil}
\newtoks\baspagetitre       \baspagetitre={\hfil}

\newtoks\hautpagegauche
\hautpagegauche={\eightpoint\rlap{\folio}\hfil\the\auteurcourant\hfil}
\newtoks\hautpagedroite
\hautpagedroite={\eightpoint\hfil\the\titrecourant\hfil\llap{\folio}}

\newtoks\baspagegauche      \baspagegauche={\hfil}
\newtoks\baspagedroite      \baspagedroite={\hfil}

\newif\ifpagetitre          \pagetitretrue


\headline={\ifpagetitre\the\hautpagetitre
\else\ifodd\pageno\the\hautpagedroite\else\the\hautpagegauche\fi\fi}

\footline={\ifpagetitre\the\baspagetitre\else
\ifodd\pageno\the\baspagedroite\else\the\baspagegauche\fi\fi
\global\pagetitrefalse}


\def\raggedbottom{\topskip 10pt plus 36pt\r@ggedbottomtrue}



\def\pointir{\unskip . --- \ignorespaces}


\def\Bigbreak{\vskip-\lastskip\bigbreak}
\def\Medbreak{\vskip-\lastskip\medbreak}


\def\ctexte#1\endctexte{%
  \hbox{$\vcenter{\halign{\hfill##\hfill\crcr#1\crcr}}$}}


\long\def\ctitre#1\endctitre{%
    \ifdim\lastskip<24pt\vskip-\lastskip\bigbreak\bigbreak\fi
  		\vbox{\parindent=0pt\leftskip=0pt plus 1fill
          \rightskip=\leftskip
          \parfillskip=0pt\bf#1\par}
    \bigskip\nobreak}

\long\def\section#1\endsection{%
\vskip 0pt plus 3\normalbaselineskip
\penalty-250
\vskip 0pt plus -3\normalbaselineskip
\Bigbreak
\message{[section \string: #1]}{\bf#1\unskip}\pointir}

\long\def\sectiona#1\endsection{%
\vskip 0pt plus 3\normalbaselineskip
\penalty-250
\vskip 0pt plus -3\normalbaselineskip
\Bigbreak
\message{[sectiona \string: #1]}%
{\bf#1}\medskip\nobreak}

\long\def\subsection#1\endsubsection{%
\Medbreak
{\it#1\unskip}\pointir}

\long\def\subsectiona#1\endsubsection{%
\Medbreak
{\it#1}\par\nobreak}

\def\rem#1\endrem{%
\Medbreak
{\it#1\unskip} : }

\def\remp#1\endrem{%
\Medbreak
{\pc #1\unskip}\pointir}

\def\rema#1\endrem{%
\Medbreak
{\it #1}\par\nobreak}

\def\newparwithcolon#1\endnewparwithcolon{
\Medbreak
{#1\unskip} : }

\def\newparwithpointir#1\endnewparwithpointir{
\Medbreak
{#1\unskip}\pointir}

\def\newpara#1\endnewpar{
\Medbreak
{#1\unskip}\smallskip\nobreak}


\long\def\th#1 #2\enonce#3\endth{%
   \Medbreak
   {\pc#1} {#2\unskip}\pointir{\it #3}\medskip}

\long\def\tha#1 #2\enonce#3\endth{%
   \Medbreak
   {\pc#1} {#2\unskip}\par\nobreak{\it #3}\medskip}


\long\def\Th#1 #2 #3\enonce#4\endth{%
   \Medbreak
   #1 {\pc#2} {#3\unskip}\pointir{\it #4}\medskip}

\long\def\Tha#1 #2 #3\enonce#4\endth{%
   \Medbreak
   #1 {\pc#2} #3\par\nobreak{\it #4}\medskip}


\def\decale#1{\smallbreak\hskip 28pt\llap{#1}\kern 5pt}
\def\decaledecale#1{\smallbreak\hskip 34pt\llap{#1}\kern 5pt}
\def\puce{\smallbreak\hskip 6pt{$\scriptstyle\bullet$}\kern 5pt}



\def\displaylinesno#1{\displ@y\halign{
\hbox to\displaywidth{$\@lign\hfil\displaystyle##\hfil$}&
\llap{$##$}\crcr#1\crcr}}


\def\ldisplaylinesno#1{\displ@y\halign{
\hbox to\displaywidth{$\@lign\hfil\displaystyle##\hfil$}&
\kern-\displaywidth\rlap{$##$}\tabskip\displaywidth\crcr#1\crcr}}


\def\eqalign#1{\null\,\vcenter{\openup\jot\m@th\ialign{
\strut\hfil$\displaystyle{##}$&$\displaystyle{{}##}$\hfil
&&\quad\strut\hfil$\displaystyle{##}$&$\displaystyle{{}##}$\hfil
\crcr#1\crcr}}\,}


\def\system#1{\left\{\null\,\vcenter{\openup1\jot\m@th
\ialign{\strut$##$&\hfil$##$&$##$\hfil&&
        \enskip$##$\enskip&\hfil$##$&$##$\hfil\crcr#1\crcr}}\right.}


\let\@ldmessage=\message

\def\message#1{{\def\pc{\string\pc\space}%
                \def\'{\string'}\def\`{\string`}%
                \def\^{\string^}\def\"{\string"}%
                \@ldmessage{#1}}}



\def\up#1{\raise 1ex\hbox{\smallf@nt#1}}


\def\qed{\raise -2pt\hbox{\vrule\vbox to 10pt{\hrule width 4pt
                 \vfill\hrule}\vrule}}

\def\cqfd{\unskip\penalty 500\quad\qed\medbreak}

\def\virg{\raise .4ex\hbox{,}}   


\def\build#1_#2^#3{\mathrel{
\mathop{\kern 0pt#1}\limits_{#2}^{#3}}}


\def\boxit#1#2{%
\setbox1=\hbox{\kern#1{#2}\kern#1}%
\dimen1=\ht1 \advance\dimen1 by #1 \dimen2=\dp1 \advance\dimen2 by #1
\setbox1=\hbox{\vrule height\dimen1 depth\dimen2\box1\vrule}%
\setbox1=\vbox{\hrule\box1\hrule}%
\advance\dimen1 by .6pt \ht1=\dimen1
\advance\dimen2 by .6pt \dp1=\dimen2  \box1\relax}


\catcode`\@=12

\showboxbreadth=-1  \showboxdepth=-1


\def\({{\rm (}}
\def\){{\rm )}}
\font\douzebf=cmbx10 at 12pt
\centerline{\douzebf  The Langlands lemma and the Betti numbers}
\centerline{\douzebf of stacks of $G$--bundles on a curve}
\vskip 10mm
\centerline{G. Laumon and M. Rapoport}
\vskip 15mm

Atiyah and Bott [AB] and Harder and Narasimhan [HN] have established
a formula for
the Poincar\'e series (the generating series formed using the Betti numbers)
of the stack ${\cal M}(G,\nu'_G)$ of $G$--bundles with slope $\nu'_G$ on a
Riemann surface,
which expresses it in terms of the Poincar\'e series  of the open substack
of semi--stable $G$--bundles and the similar Poincar\'e series for all
standard Levi
subgroups of $G$. This relation is a consequence of the Harder--Narasimhan
stratification
of ${\cal M}(G,\nu'_G)$. A similar relation arises in the context of period
domains
over a finite or $p$--adic field, where the Euler--Poincar\'e characteristic of
a generalized
flag variety of a reductive group is expressed in terms of the
Euler--Poincar\'e
characteristics
of the period domains associated with the various standard Levi subgroups of
$G$ (cf. [Rap]).  These relations
can be considered as recursion relations expressing the Poincar\'e series
(resp. Euler--Poincar\'e
characteristic) of the semi--stable sublocus in terms of the corresponding
quantities for
the ambiant spaces for $G$ and its Levi subgroups.

In this paper we show that the Langlands lemma from the theory of
Eisenstein series, which has become a standard tool in the development of the
Arthur--Selberg trace formula, can be used to invert the recursion relation for
the
Poincar\'e series of the open substack of semi--stable $G$--bundles. This note
is therefore of a purely combinatorial nature.

This application of the Langlands lemma has been noticed by Kottwitz
(in the context of $p$--adic period domains). Our only contribution has been to
formalize this
suggestion in a different context.
In the case of vector bundles on a curve the inversion of the
recursion formula had been obtained earlier by Zagier [Za] using different
techniques.
\vskip 2mm

The paper is organized as follows.
In Section 1 we fix our notations and recall the Langlands lemma. In section 2
we use the
lemma to prove a general inversion formula. In section 3 we explain how to
apply this
inversion formula to the theory of $G$--bundles on a curve. The special case
of vector bundles
is discussed in section 4.
\vskip 2mm

We wish to thank R. Kottwitz for many
helpful discussions. We also thank the Deutsche Forschungsgemeinschaft for its
support.
\vskip 5mm

{\bf 1. The Langlands lemma.}
\vskip 2mm
Let $G$ be a reductive algebraic group over a perfect field $k$. We fix a
minimal parabolic subgroup $P_0$ of $G$ and a Levi subgroup $M_0$.
We denote by ${\cal P}$ the set of
standard parabolic subgroups of $G$, i.e. parabolic subgroups of $G$ containing
$P_0$.

If $P\in {\cal P}$, we denote by $N_P$ its unipotent radical and by
$M_P$ the unique Levi subgroup of $P$ containing $M_0$. Moreover, we denote
by
$$
A_P={\rm Hom}_{k{\rm -gr}}({\bb G}_{{\rm m},k},Z_P)\otimes
{\bb G}_{{\rm m},k}
$$
the maximal split torus in the center $Z_P$ of $M_P$ and by
$$
A'_P={\rm Hom}({\rm Hom}_{k{\rm -gr}}(M_{P,\rm ab},{\bb G}_{{\rm m},k}),
{\bb G}_{{\rm m},k})
$$
the maximal quotient split torus of $M_{P,\rm ab}$. The composite map
$$
A_P\hookrightarrow Z_P\hookrightarrow  M_P\rightarrow\!\!\!\rightarrow
M_{P,\rm ab}\rightarrow\!\!\!\rightarrow A'_P
$$
is an isogeny. In particular, we have an injective map of free abelian groups
of
the same finite rank
$$
X_*(A_P)={\rm Hom}_{k{\rm -gr}}({\bb G}_{{\rm m},k},Z_P)
\hookrightarrow
{\rm Hom}({\rm Hom}_{k{\rm -gr}}(M_{P,\rm ab},{\bb G}_{{\rm m},k}),
{\bb Z})=X_*(A'_P)\ .
$$
Following Arthur [Ar1], for each $P\in {\cal P}$, we set
$$
{\goth a}_P={\bb R}\otimes X_*(A_P)={\bb R}\otimes X_*(A'_P)\ .
$$

If $P\subset Q$ are two standard parabolic subgroups of $G$, we have
canonical
maps
$$
A_Q\hookrightarrow A_P\hookrightarrow A'_P\rightarrow\!\!\!\rightarrow
A'_Q\ .
$$
The canonical maps
$A_Q\hookrightarrow A_P$ and $A'_P\rightarrow\!\!\!\rightarrow A'_Q $
induce a canonical embedding ${\goth a}_Q\hookrightarrow {\goth a}_P$ and a
canonical retraction ${\goth a}_P\rightarrow\!\!\!\rightarrow {\goth a}_Q$.
Hence, we have a canonical splitting
$$
{\goth a}_P={\goth a}_P^Q\oplus {\goth a}_Q\,,
$$
where ${\goth a}_P^Q$ is the kernel of the retraction.
Taking the dual real vector spaces, we get a splitting
$$
{\goth a}_P^*={\goth a}_P^{Q*}\oplus {\goth a}_Q^*\ .
$$
More generally, if $P\subset Q\subset R$ are three standard parabolic
subgroups of $G$,
we have  canonical splittings
$$
{\goth a}_P={\goth a}_P^Q\oplus {\goth a}_Q^R\oplus {\goth a}_R
$$
and
$$
{\goth a}_P^*={\goth a}_P^{Q*}\oplus {\goth a}_Q^{R*}\oplus {\goth a}_R^*\ .
$$
We shall denote by
$[\cdot]^Q$, $[\cdot]_Q^R$ and $[\cdot]_R$ the canonical projections of
${\goth a}_P$ onto
${\goth a}_P^Q$, ${\goth a}_Q^R$ and ${\goth a}_R$ respectively.

For each $P\in {\cal P}$, let $\Phi_P\subset {\goth a}_P^{G*}\subset
{\goth a}_P^*$
be the set of the non trivial characters of $A_P$ which
occur in the Lie algebra ${\goth g}$ of $G$ and let
$\Phi_P^+\subset\Phi_P$ be the set of the non trivial characters of $A_P$
which occur in the Lie algebra ${\goth n}_P$ of the unipotent radical of $P$.
It is well known that $\Phi_0=\Phi_{P_0}$ is a root system and
that $\Phi_0^+=\Phi_{P_0}^+$
is an order on $\Phi_0$. Let $\Delta_0=\Delta_{P_0}\subset\Phi_0^+$
be the set of simple
roots~; $\Delta_0$ is a basis of the real vector space ${\goth a}_{P_0}^{G*}$.
For each $\alpha\in\Phi_0$, there is a corresponding coroot
$\alpha^\vee$ and $(\alpha^\vee )_{\alpha\in\Delta_0}$ is a basis of the real
vector
space ${\goth a}_0^G={\goth a}_{P_0}^G\subset {\goth a}_{P_0}={\goth a}_0$.
For the other $P$'s in ${\cal P}$, $\Phi_P$ is not a root system in general.
Nevertheless,
following Arthur, we define $\Delta_P\subset\Phi_P^+$ as the set
of non trivial
restrictions to $A_P$ (or ${\goth a}_P$) of the simple roots in $\Delta_0$.
Then, $\Delta_P$ is a basis of the real vector space ${\goth a}_P^{G*}$
and, for
each $\alpha\in\Delta_P$, there is a corresponding ``coroot''
$\alpha^\vee\in {\goth a}_P^G$ with the property that
$(\alpha^\vee )_{\alpha\in\Delta_P}$
is a basis of the real vector space ${\goth a}_P^G$~: $\alpha$ is the
restriction to $A_P$ of a
unique $\beta\in\Delta_0$ and $\alpha^\vee$ is the projection of
$\beta^\vee$
onto ${\goth a}_P^G$.

If $P\subset Q$ are two standard parabolic subgroups of $G$, let
$\Phi_P^Q=\Phi_{P\cap M_Q}$
(resp. $\Phi_P^{Q+}=\Phi_{P\cap M_Q}^+$, resp. $\Delta_P^Q=
\Delta_{P\cap M_Q}$) be the set
of $\alpha$ in $\Phi_P$ (resp. $\Phi_P^+$, resp. $\Delta_P$) which
occur in the Lie
algebra ${\goth m}_Q$ of $M_Q$. Then, on the one hand, $\Delta_P^Q$ is
contained
in ${\goth a}_P^{Q*}\subset {\goth a}_P^{G*}$ and is
a basis of the real vector space ${\goth a}_P^{Q*}$. On the other hand,
the projection of
$(\alpha^\vee )_{\alpha\in\Delta_P^Q}$ onto ${\goth a}_P^Q$ is a basis of
the real
vector space ${\goth a}_P^Q$ and we may consider its dual basis
$(\varpi_\alpha^Q)_{\alpha\in\Delta_P^Q}\subset {\goth a}_P^{Q*}$.

If $P\subset Q\subset R$ are three standard parabolic subgroups of $G$ and if
$H\in {\goth a}_P^R$, we have
$$
\langle\alpha ,[H]^Q\rangle =\langle\alpha ,H\rangle\,,
\qquad\forall\alpha\in\Delta_P^Q\subset\Delta_P^R
$$
and
$$
\langle\varpi_\alpha^R,[H]_Q\rangle =\langle\varpi_\beta^R,H\rangle\,,
\qquad\forall\alpha\in\Delta_Q^R\,,
$$
where $\beta$ is the unique element in $\Delta_P^R$ such that
$\alpha =\beta |A_Q$.

\th LEMMA 1.1
\enonce
Let $P\subset R$ be two standard parabolic subgroups of $G$ and let
$H\in {\goth a}_P^R$.
\decale{\rm (i)} Let us assume that $\langle\alpha ,H\rangle >0$ or
$\langle\varpi_\alpha ,H\rangle >0$ for each $\alpha\in\Delta_P^R$.
Then, we have
$\langle\varpi_\alpha ,H\rangle >0$ for all $\alpha\in\Delta_P^R$.
\decale{\rm (ii)} Let us assume that $\langle\alpha ,H\rangle \leq 0$ or
$\langle\varpi_\alpha ,H\rangle \leq 0$ for each $\alpha\in\Delta_P^R$.
Then, we have
$\langle\varpi_\alpha ,H\rangle \leq 0$ for all $\alpha\in\Delta_P^R$.
\endth

\rem Proof
\endrem
See [La] 3.1.

\hfill\hfill\cqfd

If $P\subset Q$ are two standard parabolic subgroups of $G$, Arthur
has introduced two
characteristic functions on the real vector space ${\goth a}_P^Q$~:
the characteristic function $\tau_P^Q$ of the acute Weyl chamber
$$
{\goth a}_P^{Q+}=\{H\in {\goth a}_P^Q\mid \langle\alpha ,H\rangle >0\,,~
\forall\alpha\in\Delta_P^Q\}
$$
and the characteristic function $\widehat\tau_P^Q$ of the obtuse Weyl chamber
$$
{}^+{\goth a}_P^Q=\{H\in {\goth a}_P^Q\mid\langle
\varpi_\alpha^Q ,H\rangle >0\,,~
\forall\alpha\in\Delta_P^Q\}\ .
$$
It follows from lemma 1.1 (i) that ${\goth a}_P^{Q+}\subset {}^+{\goth a}_P^Q$.

\th LEMMA 1.2 (Langlands)
\enonce
For any standard
parabolic subgroups $P\subset R$ of $G$ and any $H\in {\goth a}_P^R$, we have
$$
\sum_{P\subset Q\subset R}(-1)^{{\rm dim}({\goth a}_Q^R)}\tau_P^Q([H]^Q)
\widehat\tau_Q^R([H]_Q)=\delta_P^R\ .
$$
and
$$
\sum_{P\subset Q\subset R}(-1)^{{\rm dim}({\goth a}_P^Q)}
\widehat\tau_P^Q([H]^Q)\tau_Q^R([H]_Q)=\delta_P^R\ .
$$
\endth

\rem Proof
\endrem
See [Ar1] \S 6 or [La] 3.2.

\hfill\hfill\cqfd

Following Arthur (see [Ar2] \S 2), we set
$$
\Gamma_P^R(H,T)=\sum_{P\subset Q\subset R}(-1)^{{\rm dim}({\goth a}_Q^R)}
\tau_P^Q([H]^Q)\widehat\tau_Q^R([H-T]_Q)
$$
and
$$
\widehat\Gamma_P^R(H,T)=\sum_{P\subset Q\subset R}
(-1)^{{\rm dim}({\goth a}_P^Q)}
\tau_P^Q([H-T]^Q)\widehat\tau_Q^R([H]_Q)=(-1)^{{\rm dim}({\goth a}_P^R)}
\Gamma_P^R(H-T,-T)\,,
$$
for any standard parabolic subgroups $P\subset R$ of $G$ and for any
$H,T\in {\goth a}_P^R$.
As an immediate consequence of the Langlands lemma, we obtain
$$
\sum_{P\subset Q\subset R}(-1)^{{\rm dim}({\goth a}_Q^R)}
\Gamma_P^Q(H,T)\widehat\Gamma_Q^R(H,T)=\delta_P^R
$$
and
$$
\sum_{P\subset Q\subset R}(-1)^{{\rm dim}({\goth a}_P^Q)}
\widehat\Gamma_P^Q(H,T)\Gamma_Q^R(H,T)=\delta_P^R\ .
$$

\th LEMMA 1.3 (Arthur)
\enonce
If $T\in {\goth a}_P^{R+}\subset {}^+{\goth a}_P^R$, the function
$H\mapsto\Gamma_P^R(H,T)$ {\rm (}resp. $H\mapsto
\widehat\Gamma_P^R(H,T)${\rm )}
is the characteristic function of the bounded subset
$$
\{H\in {\goth a}_P^R\mid \langle\alpha ,H\rangle >0\,,~
\langle\varpi_\alpha ,H\rangle\leq\langle\varpi_\alpha ,T\rangle \,,~
\forall\alpha\in\Delta_P^R\}\subset {\goth a}_P^{R+}
$$
{\rm (}resp.
$$
\{H\in {\goth a}_P^R\mid \langle\varpi_\alpha^R,H\rangle >0\,,~
\langle\alpha ,H\rangle\leq\langle\alpha ,T\rangle \,,~
\forall\alpha\in\Delta_P^R\}\subset {}^+{\goth a}_P^R\,)
$$
of ${\goth a}_P^R$.
\endth

In particular, when $T$ goes to infinity, the supports of the functions
$H\mapsto\Gamma_P^R(H,T)$ (resp. $H\mapsto\widehat\Gamma_P^R(H,T)$)
cover
${\goth a}_P^{R+}$ (resp. ${}^+{\goth a}_P^R$) (by definition,
$T\in {\goth a}_P^{R+}$
{\it goes to infinity} if $\langle\alpha ,T\rangle >0$ goes to infinity, for
each
$\alpha\in\Delta_P^R$).

\rem Proof
\endrem
Let us prove the statement about $\Gamma_P^R(H,T)$.
Let us fix $H$ and let us set
$$
I=\{\alpha\in\Delta_P^R\mid\langle\varpi_\alpha^R,H-T\rangle\leq 0\}
$$
and
$$
J=\{\alpha\in\Delta_P^R\mid \langle\alpha ,H\rangle >0\}\ .
$$
We clearly have
$$
\Gamma_P^R(H,T)=(-1)^{|\Delta_P^R-I|}
\sum_{\scriptstyle P\subset Q\subset R\atop\scriptstyle
I\subset\Delta_P^Q\subset J}(-1)^{|\Delta_P^Q-I|}=
(-1)^{|\Delta_P^R-I|}\delta_I^J\ .
$$
Therefore, $\Gamma_P^R(H,T)\not=0$ if and only if $I=J$.
Now, if $I=J$, we have
$$
\langle\alpha ,H-T\rangle\leq -\langle\alpha ,T\rangle <0\,,
\qquad\forall\alpha\in \Delta_P^R-I
$$
and
$$
\langle\varpi_\alpha^R,H-T\rangle\leq 0\,,\qquad\forall\alpha\in I\,,
$$
so that
$$
\langle\varpi_\alpha^R,H-T\rangle\leq 0\,,\qquad\forall\alpha\in\Delta_P^R
$$
by the lemma 1.1 (ii). Therefore, if $I=J$, we have $I=J=\Delta_P^R$
and $H$ satifies
the relations
$$
\langle\alpha ,H\rangle >0
$$
and
$$
\langle\varpi_\alpha ,H\rangle \leq\langle\varpi_\alpha ,T\rangle
$$
for all $\alpha\in\Delta_P^R$.
Conversely, if $H$ satisfies these relations, it is obvious that
$I=J=\Delta_P^R$.

The proof of the statement about
$\widehat\Gamma_P^R(H,T)$ is similar.

\hfill\hfill\cqfd
\vskip 5mm

{\bf 2. A general inversion formula.}
\vskip 2mm

We denote by ${\goth P}$ the set of pairs $(P,\nu'_P)$,
where $P\in {\cal P}$ and
$\nu'_P\in X_*(A'_P)$. We fix a topological abelian group $A$. A function
$$
a:{\goth P}\rightarrow A
$$
is said to be {\it $\widehat\Gamma$--converging} if it has the following
property~:

{\it For each standard parabolic subgroup
$P\subset Q$ of $G$ and each $\nu'_Q\in X_*(A'_Q)$, the finite sum
$$
\sum_{\scriptstyle \nu'_P\in X_*(A'_P)\atop\scriptstyle [\nu'_P]_Q=\nu'_Q}
\widehat\Gamma_P^Q([\nu'_P]^Q,T)a(P,\nu'_P)
$$
admits a limit as $T\in {\goth a}_P^{Q+}$ goes to infinity.}

If this is the case, we shall denote this limit by
$$
\sum_{\scriptstyle \nu'_P\in X_*(A'_P)\atop\scriptstyle [\nu'_P]_Q=\nu'_Q}
\widehat\tau_P^Q([\nu'_P]^Q)a(P,\nu'_P)\ .
$$

A function
$$
b:{\goth P}\rightarrow A
$$
is said to be {\it $\Gamma$--converging} if it has the following property~:

{\it For each standard parabolic subgroup
$P\subset Q$ of $G$ and each $\nu'_Q\in X_*(A'_Q)$, the finite sum
$$
\sum_{\scriptstyle \nu'_P\in X_*(A'_P)\atop\scriptstyle [\nu'_P]_Q=\nu'_Q}
\Gamma_P^Q([\nu'_P]^Q,T)b(P,\nu'_P)
$$
admits a limit as $T\in {\goth a}_P^{Q+}$ goes to infinity.}

If this is the case, we shall denote this limit by
$$
\sum_{\scriptstyle \nu'_P\in X_*(A'_P)\atop\scriptstyle [\nu'_P]_Q=\nu'_Q}
\tau_P^Q([\nu'_P]^Q)b(P,\nu'_P)\ .
$$

\th THEOREM 2.1
\enonce
For each $\widehat\Gamma$--converging function
$a:{\goth P}\rightarrow A$, there exists a
unique $\Gamma$--converging function $b:{\goth P}\rightarrow A$ such that,
for each $(Q,\nu'_Q)\in {\goth P}$, we have
$$
a(Q,\nu'_Q)=\sum_{\scriptstyle P\in {\cal P}\atop\scriptstyle P\subset Q}
\sum_{\scriptstyle \nu'_P\in X_*(A'_P)\atop\scriptstyle [\nu'_P]_Q=\nu'_Q}
\tau_P^Q([\nu'_P]^Q)b(P,\nu'_P)\ .
$$

The function $b$ is given by the following formula~:
for each $(Q,\nu'_Q)\in {\goth P}$, we have
$$
b(Q,\nu'_Q)=
\sum_{\scriptstyle P\in {\cal P}\atop
\scriptstyle P\subset Q}(-1)^{{\rm dim}({\goth a}_P^Q)}
\sum_{\scriptstyle \nu'_P\in X_*(A'_P)\atop\scriptstyle [\nu'_P]_Q=\nu'_Q}
\widehat\tau_P^Q([\nu'_P]^Q)a(P,\nu'_P)\ .
$$
\endth

\rem Proof
\endrem
This is an easy consequence of lemmas 1.2 and 1.3.

\hfill\hfill\cqfd

Let us now consider a particular case of this theorem which is relevant for
the computation of the Poincar\'e series of the stack of
semi--stable $G$--bundles on a curve.

For any standard parabolic subgroup $P$ of $G$, we fix
$n_P\in {\bb Z}_{\geq 0}$
and $\delta_0^P\in {\goth a}_0^{P*}\subset {\goth a}_0^*$. We assume that,
for any standard parabolic subgroups $P\subset Q$ of $G$, we have
$$
n_P\geq n_Q\,,
$$
$$
(\delta_0^Q-\delta_0^P)|{\goth a}_0^P=0
$$
and
$$
\langle\delta_P^Q,\alpha^\vee\rangle \in {\bb Z}_{>0}\qquad
(\forall\alpha\in\Delta_P^Q)
$$
where we have set
$$
\delta_P^Q=(\delta_0^Q-\delta_0^P)|{\goth a}_P^Q\ .
$$

We have
$$
\langle\delta_P^Q,[H]^Q\rangle =\langle\delta_P^G,H\rangle -\langle
\delta_Q^G,[H]_Q\rangle
$$
for every $H\in {\goth a}_P$.

We set
$$
m(P,\nu'_P)=n_P+\langle\delta_P^G,\nu'_P\rangle
$$
for each $(P,\nu'_P)\in {\goth P}$.

\th LEMMA 2.2
\enonce
Let $(Q,\nu'_Q)\in {\goth P}$.
\decale{\rm (i)} For each $(P,\nu'_P)\in {\goth P}$ such that
$P\subset Q$, $[\nu'_P]_Q=\nu'_Q$ and
$\widehat\tau_P^Q([\nu'_P]^Q)\not= 0$, we have
$$
m(P,\nu'_P)\geq m(Q,\nu'_Q)\ .
$$
\decale{\rm (ii)} For each positive integer $m$, there are only finitely many
$(P,\nu'_P)\in {\goth P}$ such that
$P\subset Q$, $[\nu'_P]_Q=\nu'_Q$,
$\widehat\tau_P^Q([\nu'_P]^Q)\not= 0$ and $m(P,\nu'_P)\leq m$.

\hfill\hfill\cqfd
\endth

We take $A$ to be a ${\bb Z}[[t]]$--module equipped with the
$t$--adic topology
and we assume that $A$ is complete for this topology. We consider an
arbitrary function
$$
a_0 :{\cal P}\rightarrow A
$$
and  set
$$
a(P,\nu'_P)=a_0(P)t^{m(P,\nu'_P)}\in A
$$
for any $(P,\nu'_P)\in {\goth P}$.
It follows from part (ii) of lemma 2.2 that the function $a$ is
$\widehat\Gamma$--converging.
Therefore, by theorem 2.1, there exists a unique $\Gamma$--converging
function
$$
b:{\goth P}\rightarrow A
$$
such that
$$
a(Q,\nu'_Q)=\sum_{\scriptstyle P\in {\cal P}\atop
\scriptstyle P\subset Q}\sum_{\scriptstyle \nu'_P\in X_*(A'_P)\atop
\scriptstyle [\nu'_P]_Q=\nu'_Q}
\tau_P^Q([\nu'_P]^Q)b(P,\nu'_P)\ .
$$
Moreover, the function $b$ is given by
$$
b(Q,\nu'_Q)=b_0(Q,\nu'_Q)t^{m(Q,\nu'_Q)}\qquad
(\forall (Q,\nu'_Q)\in {\goth P})\,,
$$
where
$$
b_0(Q,\nu'_Q)=\sum_{\scriptstyle P\in {\cal P}\atop\scriptstyle P\subset Q}
(-1)^{{\rm dim}({\goth a}_P^Q)}a_0(P)
\sum_{\scriptstyle \nu'_P\in X_*(A'_P)\atop\scriptstyle [\nu'_P]_Q=\nu'_Q}
\widehat\tau_P^Q([\nu'_P]^Q)t^{m(P,\nu'_P)-m(Q,\nu'_Q)}\in A
$$
for any $(Q,\nu'_Q)\in {\goth P}$ (cf. Lemma 2.2 (i)).

Let us consider the lattices
$$
\sum_{\alpha\in\Delta_P^Q}{\bb Z}\alpha^\vee\subset X_*(A'_P)
$$
and let us set
$$
\Lambda_P^Q =X_*(A'_P)\big/\sum_{\alpha\in\Delta_P^Q}{\bb Z}\alpha^\vee\ .
$$
Clearly, the projection $[\cdot]_Q:X_*(A'_P)\rightarrow {\goth a}_Q$
factors through $\Lambda_P^Q$
and, for each $\alpha\in\Delta_P^Q$, $\varpi_\alpha^Q:X_*(A'_P)
\rightarrow {\bb R}$
induces a homomorphism from $\Lambda_P^Q$ to ${\bb R}/{\bb Z}$.

\th LEMMA 2.3
\enonce
For each $(Q,\nu'_Q)\in {\goth P}$ and each standard parabolic subgroup
$P\subset Q$
of $G$, we have
$$\displaylines{
\qquad\sum_{\scriptstyle \nu'_P\in X_*(A'_P)\atop
\scriptstyle [\nu'_P]_Q=\nu'_Q}
\widehat\tau_P^Q([\nu'_P]^Q)t^{\langle\delta_P^Q,[\nu'_P]^Q\rangle}
\hfill\cr\hfill
=\Bigl(\prod_{\alpha\in\Delta_P^Q}
{1\over 1-t^{\langle\delta_P^Q,\alpha^\vee\rangle}}\Bigl)
\sum_{\scriptstyle\lambda\in\Lambda_P^Q\atop
\scriptstyle [\lambda]_Q=\nu'_Q}
t^{\sum_{\alpha\in\Delta_P^Q}\langle\delta_P^Q,\alpha^\vee\rangle
\langle\varpi_\alpha^Q (\lambda)\rangle}\,,\qquad}
$$
where, for each $\mu\in {\bb R}/{\bb Z}$, $\langle\mu\rangle\in {\bb R}$
is the unique representative of the class $\mu$ such that
$0<\langle\mu\rangle\leq 1$.
\endth

As $\delta_P^Q=\sum_{\alpha\in\Delta_P^Q}
\langle\delta_P^Q,\alpha^\vee\rangle\varpi_\alpha^Q$, we have
$$
\sum_{\alpha\in\Delta_P^Q}\langle\delta_P^Q,\alpha^\vee\rangle
\langle\varpi_\alpha^Q([\nu'_P]^Q+{\bb Z})\rangle\equiv
\langle \delta_P^Q,[\nu'_P]^Q\rangle\equiv 0\qquad ({\rm mod}~{\bb Z})
$$
for any $\nu'_P\in X_*(A'_P)$.

\rem Proof
\endrem
We have
$$
\sum_{\scriptstyle \nu'_P\in X_*(A'_P)\atop\scriptstyle [\nu'_P]_Q=\nu'_Q}
\widehat\tau_P^Q([\nu'_P]^Q)t^{\langle\delta_P^Q,[\nu'_P]^Q\rangle}
=\sum_{\scriptstyle \lambda\in\Lambda_P^Q\atop
\scriptstyle [\lambda]_Q=\nu'_Q}
t^{\langle\delta_P^Q,\dot\lambda \rangle}
\prod_{\alpha\in\Delta_P^Q}\sum_{\scriptstyle m_\alpha\in {\bb Z}\atop
\scriptstyle m_\alpha +\varpi_\alpha (\dot\lambda )>0}
t^{\langle\delta_P^Q,\alpha^\vee\rangle m_\alpha}
$$
where $\dot\lambda\in X_*(A'_P)$ is a representative of the class $\lambda$.
But, for each $p\in {\bb Z}_{>0}$ and each $x\in {\bb R}$, we have
$$
\sum_{\scriptstyle m\in {\bb Z}\atop\scriptstyle m+x>0}t^{pm}=
{t^{p(\langle x+{\bb Z}\rangle -x)}\over 1-t^p}\ .
$$
Hence, we have
$$
\sum_{\scriptstyle \nu'_P\in X_*(A'_P)\atop\scriptstyle [\nu'_P]_Q=\nu'_Q}
\widehat\tau_P^Q([\nu'_P]^Q)t^{\langle\delta_P^Q,[\nu'_P]^Q\rangle}
=\sum_{\scriptstyle \lambda\in\Lambda_P^Q\atop
\scriptstyle [\lambda]_Q=\nu'_Q}\prod_{\alpha\in\Delta_P^Q}
{t^{\langle\delta_P^Q,\alpha^\vee\rangle\langle \varpi_\alpha^Q
(\lambda )\rangle}
\over 1-t^{\langle\delta_P^Q,\alpha^\vee\rangle}}\ .
$$

\hfill\hfill\cqfd

To sum up, we can state~:

\th THEOREM 2.4
\enonce
There exists a unique function $b_0:{\goth P}\rightarrow A$
which satisfies the relation
$$
a_0(Q)=\sum_{\scriptstyle P\in {\cal P}\atop\scriptstyle P\subset Q}
\sum_{\scriptstyle \nu'_P\in X_*(A'_P)\atop\scriptstyle [\nu'_P]_Q=\nu'_Q}
\tau_P^Q([\nu'_P]^Q)b_0(P,\nu'_P)t^{m(P,\nu'_P)-m(Q,\nu'_Q)}\,,
$$
for each $(Q,\nu'_Q)\in {\goth P}$. This function is given by
$$\displaylines{
\quad b_0(Q,\nu'_Q)=\sum_{\scriptstyle P\in {\cal P}\atop
\scriptstyle P\subset Q}
(-1)^{{\rm dim}({\goth a}_P^Q)}a_0(P)t^{n_P-n_Q}
\Bigl(\prod_{\alpha\in\Delta_P^Q}
{1\over 1-t^{\langle\delta_P^Q,\alpha^\vee\rangle}}\Bigl)\cdot
\hfill\cr\hfill
\cdot\sum_{\scriptstyle\lambda\in\Lambda_P^Q\atop
\scriptstyle [\lambda]_Q=\nu'_Q}
t^{\sum_{\alpha\in\Delta_P^Q}\langle\delta_P^Q,\alpha^\vee\rangle
\langle\varpi_\alpha^Q (\lambda)\rangle}
\in A\,,\quad}
$$
for each $(Q,\nu'_Q)\in {\goth P}$.
\hfill\hfill\cqfd
\endth

{\pc REMARK} 2.5\pointir It follows from the definition of the function $b_0$
that $b_0(P,\nu'_P)$ only depends on
the class $\overline\nu'_P$ of $\nu'_P$ in
$ X_*(A'_P)/X_*(A_P)$.

\hfill\hfill\cqfd

\vskip 5mm

{\bf 3. Application to $G$-bundles.}
\vskip 2mm

{}From now on, let us assume that $k$ is algebraically closed,
so that $P_0$ is a Borel
subgroup of $G$, and let us fix a smooth, projective and
connected curve $X$ of genus $g\geq 2$ over $k$.

Let us recall that, for any $P\in {\cal P}$ and any $P$--bundle ${\cal T}_P$
on $X$, the {\it slope} of ${\cal T}_P$ is the element
$\mu ({\cal T}_P)\in X_*(A'_P)$ defined
by the condition
$$
\langle\xi ,\mu ({\cal T}_P)\rangle =c_1({\cal T}_{P,\xi})\in {\bb Z}\,,
\qquad\forall\xi\in X^*(A'_P)\,,
$$
where ${\cal T}_{P,\xi}$ is the line bundle on $X$ deduced
from ${\cal T}_P$ by push--out via
$P\rightarrow\!\!\!\rightarrow A'_P\buildrel\xi\over\rightarrow GL_1$.
Let us also
recall that a $P$--bundle ${\cal T}_P$ on $X$ is said to be
{\it semi--stable} (see [Ra])
if, for each standard parabolic subgroup $Q\subset P$
such that $|\Delta_Q^P|=1$
and for each $Q$--bundle ${\cal T}_Q$ on $X$ such that
$$
{\cal T}_P\cong {\cal T}_Q\times^QP\,,
$$
the slope $\mu ({\cal T}_Q)\in X_*(A'_Q)\subset {\goth a}_Q$
satisfies
$$
\langle\alpha_Q,\mu ({\cal T}_Q)\rangle\leq 0\,,
$$
where $\alpha_Q$ is the unique element of $\Delta_Q^P$.

\th LEMMA 3.1
\enonce
For each $G$--bundle ${\cal T}$ of slope $\nu'_G$, there exist
$(P,\nu'_P)\in {\goth P}$ with $[\nu'_P]_G=\nu'_G$ and
$[\nu'_P]^G\in {\goth a}_P^{G+}$,
a semi--stable $P$--bundle ${\cal T}_P$ of slope $\nu'_P$
and an isomorphism
$$
\iota :{\cal T}_P\times^PG\buildrel\sim\over\rightarrow {\cal T}\ .
$$

Moreover, the pair $(P,\nu'_P)$ and the isomorphism class of
the pair $({\cal T}_P,\iota )$
are uniquely determined by ${\cal T}$.
\endth

The pair $(P,\nu'_P)$ is called the {\it Harder--Narasimhan type}
of ${\cal T}$ and the pair
$({\cal T}_P,\iota )$ is called the {\it Harder--Narasimhan reduction}
of ${\cal T}$.

\rem Proof
\endrem
See [HN] in the case of $G=GL_n$ and [AB] in general.

\hfill\hfill\cqfd

For each $\nu'_G\in X_*(A'_G)$, we wish to consider the Poincar\'e series
$$
P_t^{\rm ss}(G,\nu'_G)\in {\bb Z}[[t]]
$$
of the stack ${\cal M}^{\rm ss}(G,\nu'_G)$ of semi--stable
$G$--bundles on $X$ of slope
$\nu'_G$. There are at least three ways to make sense of this series.
Harder and Narasimhan ([HN]) count (in a weighted way)
in the case $k=\overline {\bb F}_p$
the number of semi--stable bundles which are defined over
a finite subfield of $k$
and obtain this series as a consequence of Deligne's purity theorem.
Atiyah and Bott ([AB]) consider in the case $k={\bb C}$
the action of a gauge group on the space of complex structures
on a ${\cal C}^\infty$--bundle on the Riemann surface $X({\bb C})$ and
define $P_t^{\rm ss}(G,\nu'_G)$ as the Poincar\'e series for the
equivariant cohomology of  the semi--stable open subset.
Bifet, Ghione and Letizia ([BGL]) consider an ind-variety of semi--stable
matrix divisors and obtain  $P_t^{\rm ss}(G,\nu'_G)$ in terms of its
$\ell$--adic
cohomology.
Most probably, $P_t^{\rm ss}(G,\nu'_G)$ is also the Poincar\'e series of
the smooth algebraic stack of semi--stable $G$--bundles on $X$ of slope
$\nu'_G$ for the $\ell$--adic cohomology.

We point out that this Poincar\'e series is not the Poincar\'e polynomial
of the coarse
moduli scheme of semi--stable $G$--bundles on $X$ of slope
$\nu'_G$ (for a relation in a special case, see section 4). In fact,
it is not even a polynomial in general.

There is a recursion formula for $P_t^{\rm ss}(G,\nu'_G)$, as follows.
For each $\nu'_G\in X_*(A'_G)$, we have the stack ${\cal M}(G,\nu'_G)$ of
$G$--bundles on $X$ of slope $\nu'_G$.
For each $(P,\nu'_P)\in {\goth P}$ such that $[\nu'_P]_G=\nu'_G$ and
$[\nu'_P]^G\in {\goth a}_P^{G+}$,
we also have the substack ${\cal M}(G,P,\nu'_P)\subset
{\cal M}(G,\nu'_G)$ of $G$--bundles on $X$ of slope $\nu'_G$ which
admit $(P,\nu'_P)$
as Harder--Narasimhan type. The family of ${\cal M}(G,P,\nu'_P)$ is a
stratification
of ${\cal M}(G,\nu'_G)$, with
${\cal M}^{\rm ss}(G,\nu'_G)={\cal M}(G,G,\nu'_G)$
as the open stratum. The codimension of the stratum
${\cal M}(G,P,\nu'_P)$
is equal to
$$
{\rm dim}(N_P)(g-1)+2\langle\rho_P^G,\nu'_P\rangle\,,
$$
where
$$
\rho_P^G={1\over 2}\sum_{\alpha\in\Phi_P^{G+}}\alpha\in
{\goth a}_P^{G*}\subset
{\goth a}_P^*\ .
$$
We set
$$
m(P,\nu'_P)=2{\rm dim}(N_P)(g-1)+4\langle\rho_P^G,\nu'_P\rangle\ .
$$

We have the Poincar\'e series $P_t(G,\nu'_G)$ of ${\cal M}(G,\nu'_G)$ and
also the Poincar\'e series $P_t(G,P,\nu'_P)$ of ${\cal M}(G,P,\nu'_P)$
for any Harder--Narasimhan type $(P,\nu'_P)$.

In all the above definitions we may replace $G$ by the Levi component $M_P$
of any standard parabolic subgroup $P$ of $G$.
For each Harder--Narasimhan type $(P,\nu'_P)$ we have a fibration
$$
{\cal M}(G,P,\nu'_P)\rightarrow {\cal M}^{\rm ss}(M_P,\nu'_P)
$$
given by ${\cal T}\mapsto {\cal T}_P/N_P$, where
$({\cal T}_P,\iota )$ is the
Harder--Narasimhan reduction of ${\cal T}$.

\th THEOREM 3.2 (Harder--Narasimhan ~; Atiyah--Bott)
\enonce
The stratification of ${\cal M}(G,\nu'_G)$ by the ${\cal M}(G,P,\nu'_P)$
is perfect modulo torsion, so that for the Poincar\'e series we have
$$
P_t(G,\nu'_G)=\sum_{P\in {\cal P}}
\sum_{\scriptstyle \nu'_P\in X_*(A'_P)\atop
\scriptstyle [\nu'_P]_G=\nu'_G}
\tau_P^G([\nu'_P]^G)t^{m(P,\nu'_P)}P_t(G,P,\nu'_P)\ .
$$

Moreover, for each Harder--Narasimhan type $(P,\nu'_P)$, the above
fibration is acyclic
and we have
$$
P_t(G,P,\nu'_P)=P_t^{\rm ss}(M_P,\nu'_P)\ .
$$
\endth

Again, in this theorem, we may replace $G$ by the Levi component of
any standard
parabolic subgroup of $G$.

\rem Proof
\endrem
See [HN] and [AB] Theorem 10.10.

\hfill\hfill\cqfd

The Weyl group $W_0^G$ of $A_0$ in $G$ acts on the real vector space
${\goth a}_0^G$ and, therefore, on the graded algebra
${\rm Sym}({\goth a}_0^{G*})$
of polynomials on ${\goth a}_0^G$. It is well--known that the algebra
of invariants
$$
{\rm Sym}({\goth a}_0^{G*})^{W_0^G}
$$
(with its grading) is isomorphic to an algebra of polynomials
$$
{\bb R}[I_1,\ldots ,I_{{\rm dim}({\goth a}_0^G)}]\,,
$$
where $I_1$, ..., $I_{{\rm dim}({\goth a}_0^G)}$ are algebraically
independent homogeneous
polynomials on ${\goth a}_0^G$ of degree $\geq 2$. Let us denote by
$d_1(G),\ldots ,d_{{\rm dim}({\goth a}_0^G)}(G)$ the degrees of
these homogeneous
polynomials.
Up to a permutation, the sequence
$(d_1(G),\ldots ,d_{{\rm dim}({\goth a}_0^G)}(G))$
is canonically defined. If we set
$$
{\cal W}_G(t)=\sum_{w\in W_0^G}t^{\ell (w)}=\prod_{\alpha\in\Phi_0^{G+}}
{t^{\langle\rho_0,\alpha^\vee\rangle+1}-1\over
t^{\langle\rho_0,\alpha^\vee\rangle}-1}
$$
($\ell :W_0^G\rightarrow {\bb Z}_{\geq 0}$ is the length
function), we have
$$
{\cal W}_G(t)=\prod_{i=1}^{{\rm dim}({\goth a}_0^G)}{t^{d_i(G)}-1\over t-1}\ .
$$

\th THEOREM 3.3
\enonce
For any $\nu'_G\in X_*(A'_G)$, we have
$$
P_t(G,\nu'_G)=\Bigl({(1+t)^{2g}\over 1-t^2}\Bigr)^{{\rm dim}({\goth a}_G)}
\prod_{i=1}^{{\rm dim}({\goth a}_0^G)}
{(1+t^{2d_i(G)-1})^{2g}\over (1-t^{2d_i(G)-2})(1-t^{2d_i(G)})}\ .
$$
In particular, $P_t(G,\nu'_G)$ does not depend on $\nu'_G$.
\endth

Again, in this theorem, we may replace $G$ by the Levi component of
any standard
parabolic subgroup of $G$.

\rem Proof
\endrem
See [AB] Theorem 2.15 for the case $G=GL_n$.

\hfill\hfill\cqfd

\th THEOREM 3.4
\enonce
For any $\nu'_G\in X_*(A'_G)$, the Poincar\'e series
$P_t^{\rm ss}(G,\nu'_G)\in {\bb Z}[[t]]$
is equal to the expansion of the rational function
$$\displaylines{
\quad \sum_{P\in {\cal P}}
(-1)^{{\rm dim}({\goth a}_P^G)}
\Bigl({(1+t)^{2g}\over 1-t^2}\Bigr)^{{\rm dim}({\goth a}_P)}
\Bigl(\prod_{i=1}^{{\rm dim}({\goth a}_0^P)}
{(1+t^{2d_i(M_P)-1})^{2g}\over (1-t^{2d_i(M_P)-2})
(1-t^{2d_i(M_P)})}\Bigr)\cdot
\hfill\cr\hfill
\cdot t^{2{\rm dim}(N_P)(g-1)}
\Bigl(\prod_{\alpha\in\Delta_P}
{1\over 1-t^{4\langle\rho_P,\alpha^\vee\rangle}}\Bigl)
\sum_{\scriptstyle\lambda\in\Lambda_P^G\atop
\scriptstyle [\lambda]_G=\nu'_G}
t^{4\sum_{\alpha\in\Delta_P}\langle\rho_P,\alpha^\vee\rangle
\langle\varpi_\alpha^G (\lambda)\rangle}\quad}
$$
in ${\bb Q}(t)$.
\endth

\rem Proof
\endrem
Let
$$
a_0:{\cal P}\rightarrow {\bb Q}[[t]]
$$
be the function defined by
$$
a_0(P)=\Bigl({(1+t)^{2g}\over 1-t^2}\Bigr)^{{\rm dim}({\goth a}_P)}
\prod_{i=1}^{{\rm dim}({\goth a}_0^P)}
{(1+t^{2d_i(M_P)-1})^{2g}\over (1-t^{2d_i(M_P)-2})(1-t^{2d_i(M_P)})}\ .
$$
Our function $m(P,\nu'_P)$ is of the form
$n_P+\langle\delta_P^G,\nu'_P\rangle$
with $n_P=2{\rm dim}(N_P)(g-1)$ and $\delta_P^G=4\rho_P^G$
satisfying the hypotheses
imposed in section 2. We may therefore apply theorem 2.4.
Let $b_0$ be the unique function from ${\goth P}$ to ${\bb Z}[[t]]$ which
satisfies the relation
$$
a_0(Q)=\sum_{\scriptstyle P\in {\cal P}\atop\scriptstyle P\subset Q}
\sum_{\scriptstyle \nu'_P\in X_*(A'_P)\atop\scriptstyle [\nu'_P]_Q=\nu'_Q}
\tau_P^Q([\nu'_P]^Q)b_0(P,\nu'_P)t^{m(P,\nu'_P)-m(Q,\nu'_Q)}\,,
$$
for each $(Q,\nu'_Q)\in {\goth P}$. It follows from theorems 3.2 and 3.3 that
$$
b_0(G,\nu'_G)=P_t^{\rm ss}(G,\nu'_G)
$$
and the theorem is proved.

\hfill\hfill\cqfd

{\pc REMARK} 3.5\pointir
It follows from this theorem and remark 2.5 that $P_t^{\rm ss}(G,\nu'_G)$
only depends on the class $\overline\nu'_G$ of
$\nu'_G$ in $X_*(A'_G)/X_*(A_G)$. This can be viewed directly as follows.
Let us
arbitrarily choose a line bundle ${\cal L}$
of degree $1$ on $X$. For any $\nu_G\in X_*(A_G)$, $\nu_{G*}{\cal L}$
is an $A_G$--bundle
on $X$ and the map
$$
{\cal M}^{\rm ss}(G,\nu'_G)\rightarrow {\cal M}^{\rm ss}(G,\nu'_G +\nu_G)\,,~
{\cal T}\mapsto {\cal T}\times^{A_G}\nu_{G*}{\cal L}
$$
is an isomorphism of algebraic stacks.

\vskip 5mm
{\bf 4. The case of vector bundles.}
\vskip 2mm
Let us consider the particular case $G=GL_n$. Let us take
for $P_0$ the Borel subgroup of upper triangular matrices, so that
$A_{P_0}=A'_{P_0}=(GL_1)^n$,
${\goth a}_0={\bb R}^n$ with standard coordinates $(H_1,\ldots,
H_n)$, $\Phi_0=\{H_i-H_j\mid i\not= j\}$, $\Phi_0^+=
\{H_i-H_j\mid i< j\}$ and
$\Delta_0=\{H_i-H_{i+1}\mid i=1,\ldots ,n-1\}$.
Then, the standard parabolic subgroups of $G$ are in one to one
correspondence with
the partitions of $n$.

Let $P$ be a standard parabolic of $G$ which corresponds to
the partition $(n_1,\ldots ,n_s)$.
Then we have
$$
{\goth a}_P=\{H\in {\bb R}^n\mid H_1=\cdots =H_{n_1},\ldots ,
H_{n_1+\cdots +n_{s-1}+1}=\cdots =H_n\}\,,
$$
$$
\Delta_P=\{(H_{n_1+\cdots +n_j}-H_{n_1+\cdots +n_j+1})|{\goth a}_P\mid
j=1,\ldots ,s-1\}\,,
$$
and, for any $\alpha =
(H_{n_1+\cdots +n_j}-H_{n_1+\cdots +n_j+1})|{\goth a}_P\in\Delta_P$,
$$
\alpha^\vee
=(0,\ldots ,0,{1\over n_j},\ldots ,{1\over n_j},-{1\over n_{j+1}},\ldots ,
-{1\over n_{j+1}},
0,\ldots ,0)\,,
$$
$$
\varpi_\alpha^G=\bigl(H_1+\cdots +H_{n_1+\cdots +n_j}-
{n_1+\cdots +n_j\over n}
(H_1+\cdots +H_n)\bigr)|{\goth a}_P
$$
and
$$
\langle\rho_P,\alpha^\vee\rangle ={n_j+n_{j+1}\over 2}\ .
$$
The isomorphism
$$
{\goth a}_P\buildrel\sim\over\rightarrow {\bb R}^s\,,~(H_1,\ldots, H_n)
\mapsto (h_1,\ldots ,h_s)
$$
with $h_j=H_{n_1+\cdots +n_{j-1}+1}=\cdots =H_{n_1+\cdots +n_j}$
identifies
$$
X_*(A_P)\subset X_*(A'_P)\subset {\goth a}_P
$$
with
$$
{\bb Z}^s\subset\bigoplus_{j=1}^s{1\over n_j}{\bb Z}\subset {\bb R}^s\,,
$$
and $\alpha=(H_{n_1+\cdots +n_j}-H_{n_1+\cdots +n_j+1})|{\goth a}_P\in
\Delta_P$
with
$$
h_j-h_{j+1}\,,
$$
$\alpha^\vee$ with
$$
(0,\ldots ,0,{1\over n_j},-{1\over n_{j+1}},0,\ldots ,0)
$$
and $\varpi_\alpha^G$ with
$$
n_1h_1+\cdots n_jh_j-{n_1+\cdots +n_j\over n}
(n_1h_1+\cdots +n_sh_s)\ .
$$
Moreover, the composite map
$$
{1\over n_s}{\bb Z}\hookrightarrow
\bigoplus_{j=1}^s{1\over n_j}{\bb Z}\cong X_*(A'_P)
\rightarrow\!\!\!\rightarrow\Lambda_P^G
$$
is an isomorphism and, for any
$\lambda ={m\over n_s}\in {1\over n_s}{\bb Z}\cong \Lambda_P^G$,
we have
$$
[\lambda]_G={m\over n}\in {1\over n}{\bb Z}\cong X_*(A'_G)
$$
and
$$
\varpi_\alpha^G(\lambda )=-{n_1+\cdots +n_j\over n}m\in {\bb R}/{\bb Z}\,,~
\forall\alpha =h_j-h_{j+1}\in\Delta_P\ .
$$
Therefore, we have
$$
\prod_{\alpha\in\Delta_P}{1\over 1-t^{4\langle\rho_P,\alpha^\vee\rangle}}
=\prod_{j=1}^{s-1}{1\over 1-t^{2(n_j+n_{j+1})}}
$$
and, if $\nu'_G={d\over n}\in {1\over n}{\bb Z}\cong X_*(A'_G)$, we have
$$
\sum_{\scriptstyle\lambda\in\Lambda_P^G\atop
\scriptstyle [\lambda]_G=\nu'_G}
t^{4\sum_{\alpha\in\Delta_P}\langle\rho_P,\alpha^\vee\rangle
\langle\varpi_\alpha^G (\lambda)\rangle}=
t^{2\sum_{j=1}^{s-1}(n_j+n_{j+1})\langle -{n_1+\cdots +n_j\over n}d\rangle}
$$
(it is easy to check directly that
$\sum_{j=1}^{s-1}(n_j+n_{j+1})\langle -{n_1+\cdots +n_j\over n}d\rangle
\in {\bb Z}$).

The degrees of the invariant polynomials for $W_0^G\cong {\goth S}_n$
acting on
${\goth a}_0^G\cong {\bb R}^{n-1}$ are
$$
2,3,\ldots ,n\ .
$$
Therefore, for any $\nu'_G\in X_*(A'_G)$, we have
$$
P_t(G,\nu'_G)={(1+t)^{2g}\over 1-t^2}
\prod_{i=1}^n
{(1+t^{2i+1})^{2g}\over (1-t^{2i})(1-t^{2i+2})}\ .
$$

{}From theorem 3.4, we conclude that the Poincar\'e series
$P_t^{\rm ss}(GL_n,d/n)$
of the algebraic stack of semi--stable vector
bundles of rank $n$ and degree $d$ on the curve $X$ of genus $g\geq 2$
is equal to
$$\displaylines{
\quad \sum_{s=1}^n (-1)^{s-1}\Bigl({(1+t)^{2g}\over 1-t^2}\Bigr)^s
\sum_{\scriptstyle n_1,\ldots ,n_s\geq 1\atop\scriptstyle
n_1+\cdots +n_s=n}
\Bigl(\prod_{j=1}^s\prod_{i=1}^{n_j-1}
{(1+t^{2i+1})^{2g}\over (1-t^{2i})(1-t^{2i+2})}\Bigr)\cdot
\hfill\cr\hfill
\cdot t^{2\sum_{1\leq i<j\leq s}n_in_j(g-1)}
\Bigl(\prod_{j=1}^{s-1}{1\over 1-t^{2(n_j+n_{j+1})}}\Bigl)
t^{2\sum_{j=1}^{s-1}(n_j+n_{j+1})
\langle -{n_1+\cdots +n_j\over n}d\rangle}\ .\quad}
$$

We may also consider the stack ${\cal M}^{\rm s}(GL_n,d/n)$ of
stable vector bundles
of rank $n$ and degree $d$ on the curve $X$ (of genus $g\geq 2$).
It is an open substack
of ${\cal M}^{\rm ss}(GL_n,d/n)$, which is almost a smooth
quasi--projective variety
over $k$. More precisely, there exists a smooth
quasi--projective variety $M^{\rm s}(GL_n,d/n)$ of
dimension $(n^2-1)(g-1)$ over $k$
and a morphism of stacks ${\cal M}^{\rm s}(GL_n,d/n)\rightarrow
M^{\rm s}(GL_n,d/n)$ which is a gerb with fibers all isomorphic to $BGL_1$.

If $d$ is prime to $n$, we have
${\cal M}^{\rm s}(GL_n,d/n)={\cal M}^{\rm ss}(GL_n,d/n)$ and
$M^{\rm s}(GL_n,d/n)$ is projective over $k$. Let us denote by
$Q_t^{\rm s}(GL_n,d/n)$ the Poincar\'e polynomial of $M^{\rm s}(GL_n,d/n)$ in
this case. We have
$$
Q_t^{\rm s}(GL_n,d/n)=(1-t^2)P_t^{\rm ss}(GL_n,d/n)\ .
$$

Therefore, we have proved~:

\th THEOREM 4.1
\enonce
For each integer $d$ prime to $n$, the Poincar\'e polynomial
$Q_t^{\rm s}(GL_n,d/n)$
of the moduli space $M^{\rm s}(GL_n,d/n)$ of stable vector bundles
of rank $n$ and degree $d$ on the curve $X$ \(of genus $g\geq 2$\)
is equal to
$$\displaylines{
\qquad \sum_{s=1}^n (-1)^{s-1}{(1+t)^{2gs}\over (1-t^2)^{s-1}}
\sum_{\scriptstyle n_1,\ldots ,n_s\geq 1\atop\scriptstyle
n_1+\cdots +n_s=n}
\Bigl(\prod_{j=1}^s\prod_{i=1}^{n_j-1}
{(1+t^{2i+1})^{2g}\over (1-t^{2i})(1-t^{2i+2})}\Bigr)\cdot
\hfill\cr\hfill
\cdot t^{2\sum_{1\leq i<j\leq s}n_in_j(g-1)}
\Bigl(\prod_{j=1}^{s-1}{1\over 1-t^{2(n_j+n_{j+1})}}\Bigl)
t^{2\sum_{j=1}^{s-1}(n_j+n_{j+1})
\langle -{n_1+\cdots +n_j\over n}d\rangle}\ .\qquad}
$$

\hfill\hfill\cqfd
\endth

This last formula is equivalent to the following expression
for the Poincar\'e polynomial of
the moduli space of stable vector bundles of
rank $n$ having as determinant a fixed line bundle of
degree $d$, prime to $n$, on the curve $X$ (of genus $g\geq 2$)
$$\displaylines{
\quad \sum_{s=1}^n
(-1)^{s-1}\Bigl({(1+t)^{2g}\over 1-t^2}\Bigr)^{s-1}
\sum_{\scriptstyle n_1,\ldots ,n_s\geq 1\atop\scriptstyle
n_1+\cdots +n_s=n}
\Bigl(\prod_{j=1}^s\prod_{i=1}^{n_j-1}
{(1+t^{2i+1})^{2g}\over (1-t^{2i})(1-t^{2i+2})}\Bigr)\cdot
\hfill\cr\hfill
\cdot t^{2\sum_{1\leq i<j\leq s}n_in_j(g-1)}
\Bigl(\prod_{j=1}^{s-1}{1\over 1-t^{2(n_j+n_{j+1})}}\Bigl)
t^{2\sum_{j=1}^{s-1}(n_j+n_{j+1})
\langle -{n_1+\cdots +n_j\over n}d\rangle}\ .\quad}
$$

This was proved earlier by Zagier
using different arguments (see [Za]).
As he has remarked in loc. cit., it is not at all clear that the right
hand sides of the
last two formulas are polynomials.

Let us also point out that, if $n\geq 2$, the right hand side of
the last formula vanishes
at $t=-1$ (the order of vanishing at $t=-1$ of each summand of
the double sum is
$$
(2g-1)(s-1)+\sum_{j=1}^s\sum_{i=1}^{n_j-1}(2g-2)+
\sum_{i=1}^{s-1}(-1)=(n-1)(2g-2)\,)\ .
$$
This gives a new proof of the following result of
Narasimhan and Ramanan (see [NR])~:

\th COROLLARY 4.2 (Narasimhan and Ramanan)
\enonce
The Euler--Poincar\'e characteristic of the moduli space
of stable vector bundles of
rank $n\geq 2$ having as determinant a fixed line bundle of
degree $d$ prime to $n$ on the curve $X$ \(of genus $g\geq 2$\) is
equal to $0$.

\hfill\hfill\cqfd
\endth
\vskip 10mm

\parindent=10mm

\newtoks\ref
\newtoks\auteur
\newtoks\titre
\newtoks\editeur
\newtoks\annee
\newtoks\revue
\newtoks\tome
\newtoks\pages
\newtoks\reste
\newtoks\autre

\def\livre{\leavevmode
\llap{[\the\ref]\enspace}%
\the\auteur\pointir
{\sl \the\titre},
\the\editeur,
{\the\annee}.
\smallskip
\filbreak}

\def\article{\leavevmode
\llap{[\the\ref]\enspace}%
\the\auteur\pointir
\the\titre,
{\sl\the\revue}
{\bf\the\tome},
({\the\annee}),
\the\pages.
\smallskip
\filbreak}

\def\autre{\leavevmode
\llap{[\the\ref]\enspace}%
\the\auteur\pointir
\the\reste.
\smallskip
\filbreak}


\ref={Ar1}
\auteur={J. {\pc ARTHUR}}
\titre={A trace formula for reductive groups I: Terms associated
to classes in $G({\bb Q})$}
\revue={Duke Math. J.}
\tome={45}
\annee={1978}
\pages={911-953}
\article

\ref={Ar2}
\auteur={J. {\pc ARTHUR}}
\titre={The trace formula in invariant form}
\revue={Ann. of math.}
\tome={114}
\annee={1981}
\pages={1-74}
\article

\ref={AB}
\auteur={M.F. {\pc ATIYAH} and R. {\pc BOTT}}
\titre={The Yang-Mills equations over Riemann surfaces}
\revue={Phil. Trans. Roy. Soc. London}
\tome={A 308}
\annee={1982}
\pages={523-615}
\article

\ref={BGL}
\auteur={E. {\pc BIFET}, F. {\pc GHIONE} and M. {\pc LETIZIA}}
\titre={On the Abel--Jacobi map for divisors of higher rank on a curve}
\revue={Math. Ann.}
\tome={299}
\annee={1994}
\pages={641-672}
\article

\ref={HN}
\auteur={G. {\pc HARDER} and M.S. {\pc NARASIMHAN}}
\titre={On the cohomology groups of moduli spaces of vector bundles on
curves}
\revue={Math. Ann.}
\tome={212}
\annee={1975}
\pages={215-248}
\article

\ref={La}
\auteur={J.-P. {\pc LABESSE}}
\titre={La formule des traces d'Arthur-Selberg, S\'eminaire Bourbaki,
1984-85}
\revue={Ast\'erisque}
\tome={133-134}
\annee={1986}
\pages={73-88}
\article

\ref={NR}
\auteur={M.S. {\pc NARASIMHAN} and S. {\pc RAMANAN}}
\titre={Generalized Prym varieties as fixed points}
\revue={J. of the Indian Math. Soc.}
\tome={39}
\annee={1975}
\pages={1-19}
\article

\ref={Ra}
\auteur={A. {\pc RAMANATHAN}}
\titre={Stable principal bundles on a compact Riemann surface}
\revue={Math. Ann.}
\tome={213}
\annee={1975}
\pages={129-152}
\article

\ref={Rap}
\auteur={M. {\pc RAPOPORT}}
\reste={Non--archimedean period domains,
{\it Proceedings ICM} 1994, Z\"urich, to appear}
\autre

\ref={Za}
\auteur={D. {\pc ZAGIER}}
\titre={Elementary aspects of the Verlinde formula and of the
Harder-Narasimhan-Atiyah-Bott formula}
\revue={Israel Mathematical Conference Proceedings}
\tome={}
\annee={1995}
\pages={to appear}
\article

\vskip 10mm
\let\+=\tabalign

\line{\hbox{\kern 5mm{\vtop{\+ G\'erard LAUMON\cr
\+ URA D752 du CNRS\cr
\+ Universit\'e de Paris-Sud\cr
\+ Math\'ematiques, b\^at. 425\cr
\+ 91405 ORSAY C\'edex (France)\cr}}}
\hfill\hbox{{\vtop{\+ Michael RAPOPORT\cr
\+ FB Mathematik\cr
\+ Bergische Universit\"at\cr
\+ Gaussstrasse 20\cr
\+ 42097 WUPPERTAL (Germany)\cr}}\kern 5mm}}

\bye